\documentclass [12pt] {report}
\newcommand {\D}[2] {\displaystyle\frac{\partial{#1}}{\partial{#2}}}

\newcommand {\de} {\delta}

\newcommand {\fr} {\displaystyle\frac}
\newcommand {\wt} {\widetilde}
\newcommand {\be} {\begin{equation}}
\newcommand {\ee} {\end{equation}}
\newcommand {\ba} {\begin{array}}
\newcommand {\ea} {\end{array}}
\newcommand {\bp} {\begin{picture}}
\newcommand {\ep} {\end{picture}}
\newcommand {\bc} {\begin{center}}
\newcommand {\ec} {\end{center}}
\newcommand {\bt} {\begin{tabular}}
\newcommand {\et} {\end{tabular}}
\newcommand {\lf} {\left}
\newcommand {\rg} {\right}

\newcommand {\cR} {{\cal R}}

\newcommand {\ses} {\medskip}

\newcommand {\e} {\mathop{\rm e}\nolimits}

\newcommand {\arccosh} {\mathop{\rm arccosh}\nolimits}

\newcommand {\bibit} {\bibitem}
\newcommand {\nin} {\noindent}

 \usepackage{amsmath}

\def\2#1#2#3{{#1}_{#2}\hspace{0pt}^{#3}}
\def\3#1#2#3#4{{#1}_{#2}\hspace{0pt}^{#3}\hspace{0pt}_{#4}}
\newcounter{sctn}
\def\sec#1.#2\par{\setcounter{sctn}{#1}\setcounter{equation}{0}
                  \noindent{\bf\boldmath#1.#2}\bigskip\par}

    \newcommand {\s} {\ses}

\newcommand {\ff} {\vspace{1cm}}

\begin {document}

\begin {titlepage}

\vspace{0.1in}

\begin{center}

{\Large \bf  Pseudo-Finsleroid metric function of spatially anisotropic relativistic type}

\end{center}

\vspace{0.3in}

\begin{center}

\vspace{.15in} {\large G.S. Asanov\\} \vspace{.25in}
{\it Division of Theoretical Physics, Moscow State University\\
119992 Moscow, Russia\\
{\rm (}e-mail: asanov@newmail.ru{\rm )}} \vspace{.05in}

\end{center}

\begin{abstract}

\ses

The  paper  contributes  to  the important and urgent  problem  to extend the physical theory of  space-time in a Finsler-type way  under  the  assumption   that  the  isotropy  of  space  is  violated by   a  single  geometrically distinguished spatial direction which  destroys  the pseudo-Euclidean geometric nature of the relativistic metric and space. It proves possible to retain the fundamental geometrical property that the indicatrix should be of  the constant  curvature. Similar property appears to hold in the three-dimensional section space. The last property was the characteristic of three-dimensional  positive-definite Finsleroid space proposed and  developed in the previous work, so that  the present paper  lifts that space to the four-dimensional relativistic level. The respective  pseudo-Finsleroid   metric  function is indicated. Numerous  significant tensorial and geometrical consequences have been elucidated.

\ses

\end{abstract}

\end{titlepage}

\vskip 1cm

\setcounter{sctn}{1}

\setcounter{equation}{0}

\vspace{0.3in}

\nin
{\bf 1. Introduction}

\ff

It is well-known  that  the relativistic pseudo-Euclidean metric  possesses the geometrical  property that  the  involved  indicatrix  (mass-shell) is of the constant unit  curvature (equal to $-1$), in short $\cR_{pseudo-Euclidean ~ indicatrix}= - 1.$ This  important  property  gives  rise to numerous  important  relativistic  implications.
They can attractively be extended if the Finsler metric function
entails  the constant   curvature of indicatrix, in which case we call the function the {\it Finsleroid} metric function and denote by $F$. Respectively, in the present paper  we are aimed to have
\be
\cR_{ pseudo-Finsleroid ~ indicatrix }= - H^2.
\ee
In each tangent space, the  $H$   is a constant.

The entailed Finslerian metric tensor   is implied to be of the time-space  signature: $(+ - -\dots)$. The approach is $N$-dimensional;  N=4 in the proper space-time context.

We shall specify  the space anisotropy by  the attractive condition  that the rotational symmetry of the space is really violated  by presence of a single {\it geometrically distinguished} vector  of spacelike nature, that is to say, the anisotropy bears
the {\it axial} character.

Presence of such  anisotropy   is typical of many pictures appeared in physical applications.
  Can an occurred preferred spatial vector field, to be denoted  as $ i_{\{3\}i}(x)$,   be geometrically distinguished, in the sense that the vector $i_{\{3\}i}(x)$ can  make traces on dependence of the fundamental metric function $F(x,y)$ on local directions at the point $x$?
If resolved in positive, the deviation of the $F$ from the pseudo-Euclidean metric will influence everywhere the relativistic kinematic and dynamic, as well as the equations for the relativistic physical fields.

 The Finsler geometry [1,2] (against the Riemannian geometry [3])  presents various  genius and systematic methods to incorporate in the theory any possible anisotropy of space-time.
However, to apply the methods we must specify the Finsler-type metric function $F$ which can adequately reflect and underline the anisotropy occurred.

In the previous work [4-8], we introduced and developed the spatially anisotropic three-dimensional Finsler metric function which fulfils the properties of positive-definiteness and constancy of the indicatrix curvature. Numerous remarkable and outstanding properties of the arisen space $F^{PD}_3$ has been studied. However, the possibility to extend the space to a four-dimensional pseudo-Finsleroid relativistic space $F^{relativistic}_4$
was remained a hidden matter.
 The principle difficulty to proceed was rooted in the circumstance that, despite the Finsleroid metric function $F_3$ of the space   $F^{PD}_3$ is constructed from tangent vectors in a pure algebraic (and rather simple) way,  no four-dimensional relativistic pseudo-Finsleroid metric function $F_4$ of the algebraic form can be found (and probably does not exist) to continue  $F_3$ in  a four-dimensional framework.

The four-dimensional pseudo-Riemannian space involves the two-fold indicatrix-constancy  property. Namely, the constancy of curvature of the indicatrix (three-dimensional pseudo-sphere) as well as the spatial section of the indicatrix. The section is the two-dimensional sphere. Surprisingly, it proves possible to retain the two-fold indicatrix-constancy  property
under extending $F_3$ to a four-dimensional level. The respective pseudo-Finsleroid metric function is proposed and described in  the present paper.

A lucky and attractive possibility occurs when we undertake the systematic
analysis of the possibility to reveal a respective function  $F_4$ in terms of the
involved angles $\eta, \theta,\phi$. They play the role of the hyperbolic, azimuthal,
and polar angle, respectively, and extend the  angle triple used frequently in the pseudo-Euclidean relativistic space. We  restrict the anisotropy to the uni-directional
spatial type, such that at each point of the underground manifold exists a preferred vector of
spatial type. We naturally use the direction of the vector to be the polar axis.

In Section 2,  the initial  notions and underlying concepts  are explained.

In Section 3, the required Finsleroid  metric function $F=F(x,y)$ is
constructed by the help of the angle set $\eta, \theta,\phi$.
The polar angle  $\phi$ is given by the function $\phi=\phi(x,y)$  of simple traditional form. The explicit function  $\theta=\theta(x,y)$ can be found for the azimuthal angle $\theta$  from the function
$\tan(x,\theta)$ written in Eq. (3.20). The functions  $\eta=\eta(x,y)$   and $\theta=\theta(x,y)$ extend  their pseudo-Euclidean precursors by involving the characteristic scalars $p$ and $H$.
At the same time, the dependence  $\eta=\eta(x,y)$ cannot be written in an explicit algebraic form. The reason is that  the  explicit inverse representations  are admitted by all the structural functions entered the Finsleroid metric   except for one  function. The algebraic form of that  function ($r= r(x,\eta)$ presented in Section 3) is not sufficiently simple to perform an explicit inversion.

 This peculiarity refers also to the Finsleroid  metric function $F=F(x,y)$, since in the basic relativistic representation  $F=y^0V$ the $V$ , like to the case of the pseudo-Euclidean relativistic framework, is a function of the relativistic angle $\eta$. The very dependence $V=V(x,\eta)$ is found explicitly upon evaluation of a simple integral.

In the end of Section 3, all components of the angle representation
 $$l^i=l^i(x;\eta, \theta,\phi) $$
  of  the unit vector $l^i=y^i/F$  will be  presented in a simple algebraic form.

Nevertheless, the derivatives $\partial\eta/\partial y^i$, as well as $\partial\theta/\partial y^i$, and $\partial\phi/\partial y^i$,   admit simple algebraic representations in terms of the angle triple  $\eta, \theta,\phi.$ The latter property opens up a direct way to evaluate  the components of the  unit vector and Finsleroid metric tensor in a concise  form. With these key objects at hands, we prove able to evaluate the indicatrix metric tensor. The tensor, being written in the angle representation, just reveals the property that the {\it indicatrix is a space of the constant curvature of the value} $-H^2$  (versus the  pseudo-Euclidean  value -1). The property is the principle result of the present paper.

Another fundamental geometrical significance of the obtained Finsleroid metric function is that the three-dimensional section space reveals the property that the generalized unit sphere
(obtained in the horizontal section of the indicatrix)
is a  space of the
positive-definite constant curvature of the value  $p^2$  (versus the  Euclidean sphere curvature  value 1). The property was the characteristic of three-dimensional Finsleroid space  $F^{PD}_3$. In fact, the present paper {\it lifts} that space to the four-dimensional relativistic level.
 The $H$ and $p$ play the role of the characteristic parameters of the  Finsleroid extension of the ordinary pseudo-Euclidean relativistic space.

In Section 4, all the components of the entailed unit vector and angular metric tensor have been evaluated.

In Section 5, we find out the angle derivatives of the unit vectors to evaluate the indicatrix metric tensor.
The obtained result is  surprisingly simple. Namely,
the entailed line element  is  the $1/H$-factor of the conventional line element  over the pseudo-Euclidean relativistic hyperboloid.
Thus we demonstrate  the remarkable phenomenon that $1/H$ plays the role of the  radius of the
pseudo-Finsleroid  relativistic hyperboloid, whence the  pseudo-Finsleroid indicatrix is indeed
the space of the constant curvature value $-H^2$.
Finally, we write down the result of
 evaluation of the determinant of the respective pseudo-Finsleroid metric tensor
$g_{ij}$.

In Conclusions we emphasize several   important aspects of the developed theory.

In the case $p=1$,
the proposed function $F$ does not involve any spatial directions and, therefore, is spatially isotropic.
In Appendix
we present the evaluation chain which optimally leads to the conclusion  that in this case
the function $F$ is precisely the pseudo-Finsleroid relativistic metric function introduced and applied
in the previous work [5-7,9,10].

\ff

\setcounter{sctn}{2}            \setcounter{equation}{0}

\nin    {\bf 2.  Preliminaries}

\ff

Let $M$ be an $N$-dimensional
$C^{\infty}$
differentiable  manifold, $ T_xM$ denote the tangent space to $M$ at a point $x\in M$,
and $y\in T_xM$  mean tangent vectors.
Suppose we are given on $M$ a pseudo-Riemannian metric ${\cal S}=S(x,y)$.
 Denote by
$\cR_N=(M,{\cal S})$
the obtained $N$-dimensional pseudo-Riemannian space.
Let us also assume that the manifold $M$ admits a non-vanishing time-like 1-form
$ b= b(x,y)$
and a non-vanishing  1-form
$ i_{\{3\}}= i_{\{3\}}(x,y)$
of the space-like type.
The pseudo-Riemannian norms
are taken to be unit:
$
|| b||=1, \qquad ||i_{\{3\}}||=-1.
$
Relative  to  natural local coordinates in the space
$\cR_N$
we have the local representations
$
 b=b_i(x)y^i, \qquad   i_{\{3\}}=i_{\{3\}i}(x)y^i, \qquad
S^2= a_{ij}(x)y^iy^j,
$
with a pseudo-Riemannian metric tensor $a_{ij}(x)$, so that
the {\it time-space signature}
$
{\rm sign}(a_{ij})=(+ - -\dots)
$
takes place.

 We also  introduce the {\it transversal tensor}
  $p_{ij}=p_{ij}(x)$ as follows:
\be
p_{ij}~:=-a_{ij}+ b_ib_j -i_{\{3\}i} i_{\{3\}j},
\ee
and construct the quadratic form  $p_{ij}y^iy^j$
in terms of which we introduce the variable
\be
y_{\perp}=\sqrt{p_{ij}y^iy^j},
\ee
obtaining the decomposition
\be
S^2=    b^2-\lf(i_{\{3\}}\rg)^2 -(y_{\perp})^2.
\ee

Below, we confine our treatment to the four-dimensional case, $N=4$.
The transversal tensor can conveniently by spanned by a vector pair $(i_i, j_i)$:
\be
p_{ij} =i_ii_j +j_ij_j .
\ee
The pseudo-Riemannian metric tensor is decomposing to read
\be
a_{ij}=b_ib_j-i_ii_j-j_ij_j- i_{\{3\}i} i_{\{3\}j}.
\ee

The contravariant  reciprocal tensor
\be
a^{ij}=b^ib^j-i^ii^j-j^ij^j- i_{\{3\}}^i i_{\{3\}}^j
\ee
obeys  the reciprocity
$
a^{ij}a_{jn}=\de^i_n,
$
 where $\de^i{}_n$ stands for the Kronecker symbol.
The unit norm condition reads
\be
a^{ij}b_ib_j=1, \quad  a^{ij}i_ii_j=-1, \quad  a^{ij}j_ij_j=-1,
 \quad
 a^{ij}i_{\{3\}i} i_{\{3\}j}=-1.
\ee
The covariant index of the vector $b_i$  will be raised
 by means of the tensorial rule
$ b^i=a^{ij}b_j,$ which inverse reads $ b_i=a_{ij}b^j.$

It is convenient to  use the 1-forms
$i=i_iy^i$ and $j=j_iy^i$, in addition to $b$ and $ i_{\{3\}}$,
obtaining the variable set
\be
 z=w_3 = \fr {i_{\{3\}}}b,  \qquad  w_1=\fr ib, \quad w_2=\fr jb,
\ee
and
\be
 t=\fr{w_1}{w_2} \equiv  \fr{c_1}{c_2}  \equiv  \fr ij, \qquad
c_1=\fr {w_1}z, \qquad c_2=\fr {w_2}z.
\ee

Also, we shall apply the notation

\be
w_{\perp}=\sqrt{(w_1)^2+(w_2)^2}.
\ee

At each point $x\in M$, it proves convenient
to use the triaple
$\eta, \theta, \phi$
with the meaning of the relativistic,  azimuthal, and polar angle, respectively.
It is the direction of the vector  $i^i_{\{3\}}(x)$
that is regarded as the polar axis at point $x$.

\ff

\setcounter{sctn}{3} \setcounter{equation}{0}

\nin
{\bf 3. The pseudo-Finsleroid metric function to geometrize space-time}

\ff

Our guiding idea is to  construct the desired Finsleroid metric function
\be
F=bV(x,r)
\ee
by the help of the variable
\be
r=r(x,w_3,w_{\perp})
\ee
to be specified  in the following {\it axial} way:
\be
r=zU(x,f), \qquad f=c_2Z(x,t),
\ee
subject to the angle representation
\be
r=r(x,\eta), \qquad f=f(x,\theta),  \quad Z=Z(x,\phi),
\ee
specified by the following differential equations:
\be
\fr1V
V_{\eta}=   -   \fr1{H^2}   \fr{  \sinh\eta}   {R_1 },  \qquad
\fr 1U U_{\theta}=
\fr{1}{p^2}
\fr{   \sin\theta}
{
R_2},
\ee
were
\be
 R_1= \cosh\eta  +\sqrt{1-\fr1{p^2}+\lf(1-\fr1{H^2}\rg)  \sinh^2\eta}
\ee
and
\be
R_2= \cos\theta  +  \sqrt { \fr{1}{p^2}-1}   \sin\theta .
\ee
We assume
$$
|p|\le1, \qquad |H|\ge1.
$$

\s

The equations can explicitly be integrated:
\be
V=
C_1(x)
\fr1
{R_1}
J
\ee
\ses\\
with
\be
J=
\exp\lf(
\fr
{\sqrt{H^2-1}}
{H}
\arccosh \lf( p\fr{\sqrt{H^2-1}}{\sqrt{H^2-p^2}}\,\cosh\eta\rg)
\rg),
\ee
\ses\\
or
\be
J=
 \lf(
\sqrt{1-\fr1{H^2}}
 \,\cosh\eta
+
\sqrt{1-\fr1{p^2}+\lf(1-\fr1{H^2}\rg)  \sinh^2\eta}
\rg)^{\sqrt{1-\fr1{H^2}}},
\ee
and also
\be
U=
C_2(x)
\fr1
{
R_2}
I
\ee
\ses\\
with
\be
I=
\exp\lf(
\sqrt { \fr{1}{p^2}-1}\,\,\theta\rg).
\ee

Next, we set forth the representation
\be
\ln r =
\fr1{p^2}
\int
\fr{  d\eta}
{
R_1
\sinh\eta
},
\ee
which entails
\be
\fr1r r_{\eta}
=
\fr1{p^2}
\fr{ 1}
{
R_1\sinh\eta
}, \qquad
\eta_r=
\fr{p^2}r
R_1\sinh\eta.
\ee

\s

By making the substitution
\be
 r=
C_2(x)
\fr {\sinh\eta}{R_1}
\,
Y_1,
\ee
we can conclude that  the function $Y_1$ must obey the equation
\be
(\ln  Y_1)_{\eta}=
\lf( \fr1{p^2}-1\rg)
\fr1{\sinh\eta
  \sqrt{1-\fr1{p^2}+\lf(1-\fr1{H^2}\rg)\sinh^2\eta}}.
\ee
\ses\\
The solution reads
\be
 Y_1=
\exp\lf(
-\fr12
\sqrt{ \fr1{p^2}-1}
\arctan \lf(\sqrt { \fr1{p^2}-1}
\fr{2\cosh\eta\sqrt{\fr1{H^2}-\fr1{p^2}+\lf(1-\fr1{H^2}\rg)\cosh^2\eta}}
{\fr1{H^2}-\fr1{p^2}+\lf(\fr1{H^2}-\fr1{p^2}+2\lf(1-\fr1{H^2}\rg)\rg)\cosh^2\eta
}
\rg)\rg).
\ee

After that, we need to apply the integral
\be
\ln f
=
\int
\fr{ d\theta}
{
\sin\theta
\lf(
\cos\theta
+
\sqrt { \fr{1}{p^2}-1}
 \sin\theta
\rg)
},
\ee
so that
$$
\D{\theta}f
=
\fr1f
\sin\theta
\lf(
\cos\theta
+
\sqrt { \fr{1}{p^2}-1}
 \sin\theta
\rg),
\qquad
\theta_f
=
\fr1f
\sin\theta
\fr1U
I.
$$
It can readily be verified that
\be
f
=
C_{17}(x)
\fr {\sin\theta}{\cos\theta
+
\sqrt { \fr{1}{p^2}-1}\,\, \sin\theta
}
\equiv
C_{17}(x)
\fr {\sin\theta}
{R_2}.
\ee

The last function can be inverted, yielding
\be
\tan\theta=
\fr{
\fr1{C_{17}}f
}{1-\sqrt { \fr{1}{p^2}-1} \fr1{C_{17}}f}.
\ee
The property $f(x,0)=0$ holds.

Noting also (3.11), we obtain
\be
\sin\theta=
\fr1U
f
I, \qquad
\theta_f
=
\fr1{U^2}
I^2.
\ee

\s

The polar angle $\phi$ is constructed in the traditional way:

\be
\phi=\arctan t,
\ee
so that
$$
\phi_t=\fr 1{1+t^2}.
$$
The last equality entails  $ Z=\sqrt{1+t^2}\,C_{11}(x)$  and
\be
f=\fr1{w_3}\sqrt{w_1w_1+w_2w_2} \,C_{11}(x).
\ee

Henceforth, it is convenient to specify the integration constants according to
$
C_{2}=C_{17}=1.
$
Using (3.20) and (3.21), we can write
$$
U^2=
\lf(
\lf(
1-
\sqrt { \fr{1}{p^2}-1}
f\rg)^2
+f^2
\rg)
I^2
=
\lf(
1-
2
\sqrt { \fr{1}{p^2}-1}
f
+
\fr{1}{p^2}
f^2
\rg)
I^2.
$$
Making the choice
 $ C_{11}=p,$
we come to
\be
U^2=
\lf(
1-
2
\sqrt {1-p^2}
w
+
w^2
\rg)
I^2
\ee
and
\be
f=wp,
\ee
where
\be
w=\fr1{w_3}\sqrt{w_1w_1+w_2w_2}.
\ee

With the function $U$ given by (3.24), the product $r=zU$
is just the {\it three-dimensional positive-definite Finsleroid metric function}
proposed and developed in the previous work [5-8] (where the notation $h$ was used instead of $p$).
The axial structure and the positive-definiteness of the entailed metric tensor, taken in conjunction with the indicatrix curvature constancy of positive value, is the characteristic of the function.

Now, the angle representation of tangent vectors $\{y^i\}$ can  be derived.
Indeed, at any fixed background point $x\in M$ , we have
$$
w_3=w_3(\eta,\theta), \qquad  w_{\perp}=w_{\perp}(\eta,\theta),
$$
with
$$
w_3=\fr rU      , \qquad
w_{\perp}=w_3\fr1pf    =   w_3\fr1p\fr{\sin\theta}{R_2},
$$
so that
\be
w_3 =  r\fr 1IR_2, \qquad          w_{\perp}=   \fr1p r\fr1I\sin\theta.
\ee
Also,
\be
w_1=w_{\perp}\cos\phi, \qquad    w_2=w_{\perp}\sin\phi.
\ee
Thus, the sought result is given by the equalities
$$
y^1=w_1b,  \qquad  y^2=w_2b, \qquad  y^3=w_3b,
$$
together with
\be
y^0=\fr1VF, \qquad  V=V(\eta)
\ee
(here $y^0=b$).

\ff

\setcounter{sctn}{4} \setcounter{equation}{0}

\s

\nin
{\bf 4. Unit vectors and metric tensor}

\ff

To study the geometrical aspects of the pseudo-Finsleroid space,
 we need to clarify the structure of the entailed metric tensor.

Choosing a fixed tangent space, we shall represent the involved vectors
by the help of their components with  respect to the base frame
$\{b_i, i_i, j_i,i_{\{3\}i}\}$
which enters the Riemannian metric tensor $a_{ij}$ according to (2.5).

\s

It is convenient to use
the notation
$$
V_1=\D V{w_1}, \qquad V_2=\D V{w_2}, \qquad V_3=\D V{w_3}, \qquad
V_r=\D Vr,
$$
and
$$
r_1=\D r{w_1}, \qquad r_2=\D r{w_2}, \qquad r_3=\D r{w_3}.
$$
The identity
\be
r=r_1w_1+ r_2w_2+ r_3w_3
\ee
holds fine due to the homogeneity involved.

We should derive the covariant unit vector
components $ l_i=\partial  F/\partial{y^i}$
from the function
$$
F=y^0V.
$$
We
 get
\be
l_1=V_1, \quad l_2=V_2, \quad l_3=V_3,
\ee
and
$$
l_0=V-( V_1w_1+ V_2w_2+ V_3w_3)
=
V-V_r( r_1w_1+ r_2w_2+ r_3w_3),
$$
or
$$
l_0=V-V_rr
$$
(because of (4.1)).

\s

Since
$$
\fr1VV_r=
- \fr1{H^2}
\fr{p^2}r\sinh^2\eta
$$
(see (3.5) and (3.14)),
we obtain the clear and useful representation
\be
l_0=V
\lf(1+\fr{p^2}{H^2}
\sinh^2\eta
\rg).
\ee

With this preparation,  all the derivatives
$$
l_{ij}=\D{l_{i}}{y^j}
$$
can readily be found. At first, we obtain
$$
l_{00}=\fr1{y^0}V_{rr}r^2, \qquad
l_{01}=-\fr1{y^0}V_{rr}rr_1, \qquad  l_{02}=-\fr1{y^0}V_{rr}rr_2, \qquad l_{03}=-\fr1{y^0}V_{rr}rr_3.
$$
After that, we can use
$
V_1=V_rr_1
\equiv l_1
$
and conclude that
$$
y^0l_{11}=
V_{rr}r_1r_1+V_rr_{11}.
$$

By following this method, we can find all the components of  the angular metric tensor
$$
h_{ij}=Fl_{ij}.
$$
The result reads
\be
h_{00}=VV_{rr}r^2, \quad h_{01}=-VV_{rr}rr_1, \quad  h_{02}=-VV_{rr}rr_2, \quad h_{03}=-VV_{rr}rr_3,
\ee

\s

\s

\be
h_{11}=
VV_{rr}r_1r_1+VV_rr_{11}, \qquad
h_{12}=
VV_{rr}r_1r_2+VV_rr_{12},  \qquad
h_{22}=
VV_{rr}r_2r_2+VV_rr_{22},
\ee

\s

\s

\be
h_{13}=
VV_{rr}r_1r_3+VV_rr_{13},  \qquad
h_{23}=
VV_{rr}r_2r_3+VV_rr_{23},  \qquad
h_{33}=
VV_{rr}r_3r_3+VV_rr_{33}.
\ee

Using the derivative values
$$
\fr1VV_r=
- \fr1{H^2}
\fr{p^2}r\sinh^2\eta, \quad
{\eta}_r
= \fr{p^2}r\sinh\eta
 R_1,
 \quad
\fr1{H^2}\lf({\eta}_r\rg)^2=
-\fr1VV_{rr},
 \quad
\fr1V
V_{\eta}=
-
\fr1{H^2}
\fr1{ R_1}
\sinh\eta,
$$
\ses\\
we straightforwardly obtain the angle representation
\be
h_{ij}=
-
\fr1{H^2}
\lf(
\eta_i\eta_j+\sinh^2\eta(\theta_i\theta_j+\sin^2\theta\phi_i\phi_j)
\rg)
F^2.
\ee

\ff

\setcounter{sctn}{5} \setcounter{equation}{0}

\nin
{\bf 5. The indicatrix curvature}

\ff

With the help of the formulas
written in the end of Section 3,
the differentiation of  the unit  vector components
$$
l^0=\fr1V, \quad l^a=\fr1Vw_a
$$
with respect to the used angles yields the following list:
$$
l^0_{\eta}=-\fr1VV_{\eta}l^0 ,   \quad  l^0_{\theta}= l^0_{\phi}=0,
$$

\s

$$
l^1_{\eta}=-\fr1VV_{\eta}l^1+\fr1rr_{\eta}l^1,
  \qquad
 l^1_{\theta}=\lf(-\fr1I I_{\theta}+\fr{\cos\theta}{\sin\theta}\rg)l^1,
\qquad
l^1_{\phi}=-\fr{\sin\phi}{\cos\phi}l^1,
$$

\s

$$
l^2_{\eta}=-\fr1VV_{\eta}l^2+\fr1rr_{\eta}l^2,
  \qquad
 l^2_{\theta}=\lf(-\fr1I I_{\theta}+\fr{\cos\theta}{\sin\theta}\rg)l^2,
\qquad
l^2_{\phi}=\fr{\cos\phi}{\sin\phi}l^2,
$$

\s

$$
 l^3_{\eta}=-\fr1VV_{\eta}l^3+\fr1rr_{\eta}l^3, \qquad
l^3_{\theta}= \lf(\fr1{R_2}R_{2\theta}-\fr1II_{\theta}\rg)l^3,
\qquad
l^3_{\phi}= 0,
$$
together with
$$
l^3_{\theta}= -\fr1{p^2R_2}\sin\theta l^3  = -\fr1{p^2}\sin\theta \fr1V \fr1Ir.
$$

From (3.7) and (3.12) we can conclude that
$$
\fr1II_{\theta}=\sqrt { \fr{1}{p^2}-1}, \quad
R_{2\theta}= -\sin\theta  +  \sqrt { \fr{1}{p^2}-1}   \cos\theta .
$$

\s

We straightforwardly obtain the angle representation
\be
h_{ij}=
-
\fr1{H^2}
\lf(
\eta_i\eta_j+\sinh^2\eta(\theta_i\theta_j+\sin^2\theta\phi_i\phi_j)
\rg)
F^2.
\ee

\s

Now we are prepared to evaluate  the {\it indicatrix metric tensor}
\be
i_{ab}=-h_{ij}l^i_{a}l^j_{b}, \quad (a,b=\eta,\theta, \phi)
\ee
in the angle representation.
The following  simple result  is obtained:
$$
i_{\eta\eta}= -\fr1{H^2}, \quad
i_{\theta\theta}=-\fr1{H^2} \sinh^2\eta, \quad
i_{\phi\phi}=-\fr1{H^2}\sinh^2\eta \sin^2\theta,
$$

\s

$$
i_{\eta 0}= i_{\theta 0}=_{\phi 0}=0, \qquad
i_{\eta \theta}= i_{\eta\phi}=i_{\theta\phi }=0,
$$

The entailed line element square reads
$$
(ds)^2=(dF)^2-F^2(dl)^2
$$
with
$$
(dl)^2=
\fr1{H^2}
\bigl((d \eta)^2 +\sinh^2\eta((d\theta)^2 +(d\phi)^2)\bigr).
$$
\ses\\
This is just the $1/H^2$-factor of the  metric on the pseudo-Euclidean relativistic hyperboloid.
We observe the remarkable phenomenon that $1/H$ plays the role of the {\it radius} of the
pseudo-Finsleroid  relativistic hyperboloid, whence the  pseudo-Finsleroid indicatrix is indeed
the space of the constant curvature value $-H^2$.

The attentive evaluation of the determinant of the respective pseudo-Finsleroid metric tensor
$$
g_{ij}=h_{ij}+l_il_j
$$
 leads to the following result:
\be
\det(g_{ij})
=
-\fr1{H^6}
\lf(  p^4
I^3V^4R_1
 \rg)^2
\fr1{r^6}
\sinh^6\eta.
\ee
We observe here the  independence of the  polar angle $\phi.$
The dependence on the  angle $\theta$ enters through the function $I$ written in (3.12).
The right-hand part of the determinant (5.3) tends to $-1$ when $H\to 1, p\to 1$.

\ff

\nin
{\bf Conclusions: Important aspects }

\ff

The pseudo-Finsleroid approach developed   in the present  paper  involves   two characteristic parameters of extension, $H$ and $p.$ The metric function $F$ reduces to the pseudo-Riemannian metric function when $H\to 0, p\to 0$. The  geometrical signicance of the parameters is remarkable, namely they represent  the indicatrix curvarure values. The entailed Finslerian metric tensor   is of the time-space  signature $(+---)$.
The respective  pseudo-Finsleroid   metric  function is {\it implicitly} proposed by a simple algebraic equation. However, the  components of the Finslerian unit vector and metric tensor,  as well as the partial derivatives of the angles $\eta, \theta, \phi$ are evidenced  to show sufficiently transparent and handy structure.

Unlike to the pseudo-Riemannian geometry, the  curvature  tensor  of  the  tangent  space  does not vanish  identically. The metric tensor does depend on tangent vectors: $g_{ij}=g_{ij}(x,y)$ in contrast to the Riemannian preassumption $a_{ij}=a_{ij}(x)$.

Above, we have found the function $F$ by starting from the attractive idea to preserve the key pseudo-Euclidean geometrical property that the curvature of the associated relativistic indicatrix is everywhere constant.

In the preceding papers [6-10], the pseudo-Finsleroid relativistic    metric function revealing such a property was constructed upon assuming that the involved geometrically distinguished vector (that was denoted by $b_i$) is time-like and the spatial isotropy holds fine.
In the preceding sections of the present paper it is shown that the metric admits
a due extension which also involves  the preferred vector of  space-like type.

Various  interesting physical applications of the proposed pseudo-Finsleroid metric function  can be elaborated.
Respective  extensions of the relativistic kinematic and dynamic are the nearest ways.

\ff

\setcounter{equation}{0}

\nin
{\bf Appendix: Spatially isotropic limit}

\ff

Let us put
$
p=1
$
in the representations given in Section 3. Among them, we choose the following:
$$
J=
\exp
\lf( \sqrt{1-\frac1{H^2}}  \eta  \rg), \qquad
 R_1= \cosh\eta  +\sqrt{1-\fr1{H^2}}  \sinh\eta, \quad
 r=
C_2
\fr {\sinh\eta}{R_1}.
$$
It follows that
$$
\fr1{C_2}r
\lf(\cosh\eta  +\sqrt{1-\fr1{H^2}}  \sinh\eta\rg)
=
\sinh\eta, \qquad
\fr1{C_2}r
\cosh\eta
=
\lf(
1-\fr1{C_2}r
\sqrt{1-\fr1{H^2}} \rg) \sinh\eta.
$$
This equation can be solved, yielding
$$
\fr1{\tanh\eta}
=
C_2
\fr1r
-
\sqrt{1-\fr1{H^2}}, \qquad
\fr1{\sinh^2\eta}
=
\lf(
C_2
\fr1r
-
\sqrt{1-\fr1{H^2}}
\rg)^2
-1.
$$

\s

In this way, we obtain the angle functions
$$
\sinh\eta=
\fr1
{\sqrt{
\lf(
C_2
\fr1r
-
\sqrt{1-\fr1{H^2}}
\rg)^2
-1
}
  }, \qquad
\cosh\eta=
\fr{
C_2
\fr1r
-
\sqrt{1-\fr1{H^2}}
}
{\sqrt{
\lf(
C_2
\fr1r
-
\sqrt{1-\fr1{H^2}}
\rg)^2
-1
}
  },
$$
which entail the representations
$$
R_1=
\fr{
C_2
\fr1r
}
{\sqrt{
\lf(
C_2
\fr1r
-
\sqrt{1-\fr1{H^2}}
\rg)^2
-1
}
  }, \qquad
\e^{\eta}=
\fr{
C_2
\fr1r
-
\sqrt{1-\fr1{H^2}}
+1
}
{\sqrt{
\lf(
C_2
\fr1r
-
\sqrt{1-\fr1{H^2}}
\rg)^2
-1
}
  },
$$
so that
$$
J=
\lf(
\fr{
C_2
\fr1r
-
\sqrt{1-\fr1{H^2}}
+1
}
{\sqrt{
\lf(
C_2
\fr1r
-
\sqrt{1-\fr1{H^2}}
\rg)^2
-1
}
  }
\rg)^{\sqrt{1-\fr1{H^2}}}.
$$

\ff

The function $V$ written in Section 3 reduces to
$$
V=
\fr1{ C_2  }r
\sqrt{
\lf(
C_2
\fr1r
-
\sqrt{1-\frac1{H^2}}
\rg)^2
-1
}
 \lf(
\fr{
C_2
\fr1r
-
\sqrt{1-\fr1{H^2}}
+1
}
{\sqrt{
\lf(
C_2
\fr1r
-
\sqrt{1-\fr1{H^2}}
\rg)^2
-1
}
  }
\rg)^{\sqrt{1-\fr1{H^2}}}
$$
(we have taken the integration constant $C_1$ to be 1).

\s

Let us  square this function:
$$
V^2=
\fr1{ (C_2 )^2 }r^2
\lf[
\lf(
C_2
\fr1r
-
\sqrt{1-\fr1{H^2}}
\rg)^2
-1
\rg]
 \lf(
\fr{
C_2
\fr1r
-
\sqrt{1-\fr1{H^2}}
+1
}
{
C_2
\fr1r
-
\sqrt{1-\fr1{H^2}}
-1
  }
\rg)^{\sqrt{1-\fr1{H^2}}}.
$$
\ses\\
Convenient cancelation is  possible here. We obtain
$$
V^2=
\fr1{ (C_2 )^2 }
\lf(
C_2
-
\lf(1+\sqrt{1-\fr1{H^2}}\rg)r
\rg)^{1-\sqrt{1-\fr1{H^2}}} \lf(
C_2
+\lf(1-\sqrt{1-\fr1{H^2}}\rg)r
\rg)^{1+\sqrt{1-\fr1{H^2}}}.
$$

\ff

Since
$
w_3 =  r\cos\theta,         w_{\perp}=    r\sin\theta,
$
we have
$
r=
\sqrt{(w_1)^2+(w_2)^2+(w_3)^2}\equiv |w|.
$

\s

Using the notation
\be
g_+=-\sqrt{H^2-1}+H, \qquad  g_-=-\sqrt{H^2-1}-H,
\ee
and
$\wt C_2= H C_2$,
we can write
\ses\\
\be
V^2=
\fr1{ (\wt C_2 )^2 }
\lf(
\wt C_2
+g_-|w|
\rg)^{g_+/H}
 \lf(
\wt C_2
+g_+|w|
\rg)^{-g_-/H}.
 \ee
\ses\\
This $V^2$ does coincide with the respective function used in [5-10].

\ses

\ses

\def\bibit[#1]#2\par{\rm\noindent\parskip1pt
                     \parbox[t]{.05\textwidth}{\mbox{}\hfill[#1]}\hfill
                     \parbox[t]{.925\textwidth}{\baselineskip11pt#2}\par}

\nin {  REFERENCES}

\ses

\bibit[1] H. Rund: \it The Differential Geometry of Finsler
 Spaces, \rm Springer, Berlin 1959.

\bibit[2] G.S. Asanov: \it Finsler Geometry, Relativity and Gauge
 Theories, \rm D.~Reidel Publ. Comp., Dordrecht 1985.

\bibit[3] J.L. Synge: \it Relativity: The General Theory, \rm
 North-Holland, Amsterdam I960.

\bibit[4] G.S. Asanov: {\it Moscow University Physics Bulletin:}
Class of spherically symmetric Finsler metric functions with the indicatrix of constant curvature,
 \bf34\rm(3) (1993),76-78;

\bibit[5] G.S. Asanov: Finsler cases of GF-space, \it Aeq. Math. \bf49 \rm(1995), 234-251.

\bibit[6] G.S. Asanov: Finslerian metric functions over product $R \times M$ and their potential applications,
\it Rep. Math. Phys. \bf 41\rm(1) \rm (1998), 117-132.

\bibit[7]   G. S. Asanov:
   Finsleroid--Finsler spaces of positive--definite and relativistic types,
\it Rep. Math. Phys.  \bf 58\rm(2) \rm(2006), 275-300.

\bibit[8]
G. S. Asanov:  Finsler connection preserving the two-vector  angle  under
the  indicatrix-inhomogeneous treatment,   {\it  arXiv:} 1009.2673 [math.DG] (2011);
 Finsler connection properties generated by   the two-vector angle developed on the indicatrix-inhomogeneous level,
 {\it Publ. Math. Debrecen } {\bf 82/1} (2013), 125-153.

\bibit[9] G.S. Asanov: {\it Moscow University Physics Bulletin:}
Finslerian relativistic generalization compatible with the spatial isotropy
\bf35\rm(2) (1994),19-23;
Finslerian realization of the momentum space
\bf35\rm(2) (1994),12-20;
  Finsler generalization of Lorentz transformations,
 \bf36\rm(4) (1995),7-12;
 Finslerian invariant and coordinate length in the inertial reference frame
 \bf53\rm(1) (1998),18-23.

\bibit[10]   G. S. Asanov:
   Finsleroid--Finsler spaces of positive--definite and relativistic types,
\it Rep. Math. Phys.  \bf 58\rm(2) \rm(2006), 275-300.

\end{document}